\documentclass{amsart}

\usepackage{amsmath,amssymb,amsthm}

\usepackage{tikz}
\usetikzlibrary{calc,intersections,through,backgrounds}
\usepackage[latin1]{inputenc}

\usepackage{caption}
\usepackage{graphicx}
\usepackage{graphics}

\newtheorem*{thm}{Theorem}

\begin{document}

\title[]{Fast Localization of eigenfunctions\\ via smoothed potentials}

\author[]{Jianfeng Lu}
\address[JL]{Department of Mathematics, Department of Physics, and Department of Chemistry,
Duke University, Box 90320, Durham NC 27708, USA}
\email{jianfeng@math.duke.edu}

\author[]{Cody Murphey}
\address[CM]{Department of Computer Science, Yale University, New Haven, CT 06511, USA} \email{cody.murphey@yale.edu}

\author[]{Stefan Steinerberger}
\address[SS]{Department of Mathematics, University of Washington, Seattle, WA 98195, USA} \email{steinerb@uw.edu}

\thanks{J.L.~is partially supported by the National Science Foundation via grants DMS-2012286 and CHE-2037263. S.S.~is partially supported by the NSF (DMS-1763179) and the Alfred P. Sloan Foundation.}

\begin{abstract} We study the problem of predicting highly localized low-lying eigenfunctions $(-\Delta +V) \phi = \lambda \phi$ in bounded domains $\Omega \subset \mathbb{R}^d$ for rapidly varying potentials $V$. 
 Filoche \& Mayboroda introduced the function $1/u$, where $(-\Delta + V)u=1$, as a suitable regularization of $V$ from whose minima one can predict the location of eigenfunctions with high accuracy. We proposed a fast method that produces a landscapes that is exceedingly similar, can be used for the same purposes and can be computed very efficiently: the computation time on an $n \times n$ grid, for example, is merely $\mathcal{O}(n^2 \log{n})$, the cost of two FFTs.  
\end{abstract}
\maketitle

\tableofcontents

\section{Introduction}
\subsection{The Landscape Function} Eigenfunctions of elliptic
differential operators are typically of comparable size throughout the
domain. However, underlying backgrounds composed of inhomogenious
materials will sometimes produce localized vibration patterns. Let $\Omega \in \mathbb{R}^d$ be an open, bounded domain in
which we consider
\begin{align*}
(-\Delta + V)\phi &= \lambda \phi~ \quad \mbox{in~}\Omega \\
 \phi&= 0 \qquad \mbox{on}~ \partial \Omega.
\end{align*}
Here $V:\Omega \rightarrow \mathbb{R}_{\geq 0}$ is a real-valued,
nonnegative potential which one should assume to be rapidly
oscillating or possibly changing its behavior rapidly from one region
to the next; in this setting, low-lying eigenfunctions that are
strongly localized in space become possible \cite{anderson}. Since
low-lying eigenfunctions of such a Schr\"odinger operator determine the long-time behavior of many
associated dynamical systems (say, the heat equation, the wave
equation or the Schr{\"o}dinger equation), they
are of obvious interest.
Filoche \& Mayboroda \cite{fil} provided an astonishingly effective method to predict the behavior of low-energy eigenfunctions for such operators $-\Delta + V$. They define the \textit{landscape function} as the unique function $u:\Omega: \mathbb{R} \rightarrow \mathbb{R}_{}$ solving
\begin{align*}
(-\Delta + V)u &=1~ \qquad \mbox{in~}\Omega \\
 u&= 0 \qquad  \mbox{on}~~ \partial \Omega
\end{align*}
and show that $u$ exerts pointwise control on all eigenfunctions $(-\Delta + V)\phi = \lambda \phi$ 
$$ |\phi(x)|  \leq \lambda u(x)  \|\phi\|_{L^{\infty}(\Omega)} .$$ 
An eigenfunction $\phi$ can only localize in $\left\{x: u(x) \geq 1/\lambda\right\} \subset \Omega$. This inequality, however, is not the full story: the landscape function is much more effective than is indicated by that inequality alone. Numerical experiments \cite{arnold0, arnold2, fil, fil2, fil3} suggest that the largest local maxima correspond precisely to the location where the first few eigenfunctions localize and that many more properties (including refined eigenvalue estimates and improvements on the Weyl law) are being captured. A different way of thinking about it as that the inverse $1/u$ acts as a suitably regularized potential. The accuracy of these refined predictions is quite striking and have already led to many interesting results \cite{arnold0, arnold, arnold2, chal, david, fil, fil2, fil3, har, lef, leite, pic}.

\subsection{The Universal Kernel} The main problem is to extract
the localization information out of the potential $V$. The landscape
function does this by solving $(-\Delta + V)u=1$. In a recent work
\cite{steini3}, the third author proposed an alternative approach: the
starting point is the interpretation of localized eigenstates as
critical points of
$$J(\phi) =  \int_{\Omega}{ |\nabla \phi(x)|^2 dx} + \int_{\Omega}{ V(x) \phi(x)^2 dx}$$ 
(subject to some orthogonality condition with respect to earlier
eigenfunctions). This energy functional clearly illustrates that the
relevant quantity is not the point-wise value of $V$ but, as seen in
the second integral, the integral average of $V$ over a certain
scale. One question remains: which average? If $k_t$ is a smooth,
radial, probability density centered at the origin and having most of
its mass at a ball of radius $\sim \sqrt{t}$ (the scaling is motivated
by the Theorem below), we could see what happens in our PDE when we
replace the potential $V$ by a slightly smoothed version:
$$ -\Delta \phi(x) +  (V*k_t)(x) \phi(x) = \lambda \phi(x) + \mbox{error}(x,t).$$
Certainly, the error will now depend on which type of kernel $k_t$ we choose. As it turns out, there is a unique distinguished kernel that results in the best possible dependence of the error term. 
This kernel 
 $k_t:\mathbb{R}^d \rightarrow \mathbb{R}_{\geq 0}$ is given by
 $$ k_t(x) = \frac{1}{t} \int_0^t \frac{ \exp\left( - \|x\|^2/ (4s) \right)}{(4 \pi s )^{d/2}} ds.$$
The kernel depends only on a scale parameter $t>0$ and
the dimension $d$ of the domain $\Omega \subset \mathbb{R}^d$. We note that, in practice, the boundary conditions of the domain $\Omega$ should also be relevant. However, in the setting that we consider, the potential $V$ implicitly isolates subregions: whenever $V$ is large, or `larger than average', low-frequency eigenfunctions will be small in that region, leading to boundary conditions that have a negligible effect. This is also in line with the discussion in \cite{arnold0, arnold, arnold2, fil, fil2, fil3}: the boundary conditions are only relevant near the boundary since localized solutions eigenfunctions have, unless localized near the boundary, limited interaction with the boundary.

\begin{center}
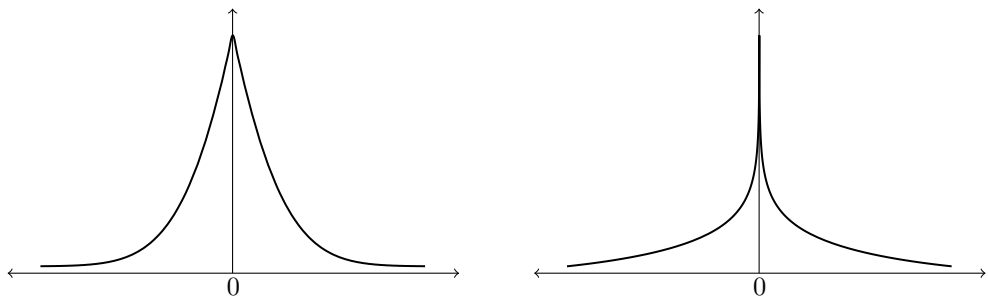
\begin{figure}[h!]
\begin{tikzpicture}
\node at (0,0.1) {\includegraphics[width=0.42\textwidth]{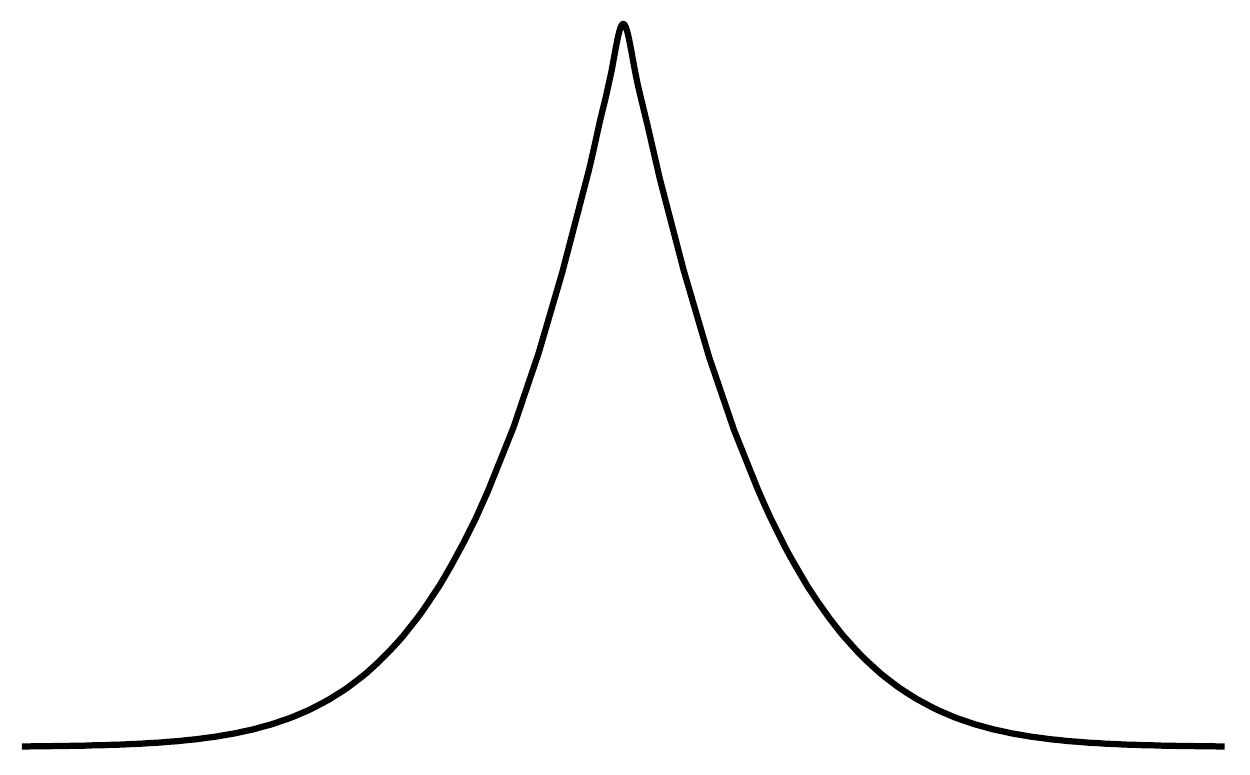}};
\draw [<->] (-3,-1.52) -- (3,-1.52);
\draw [->] (-0.01,-1.52) -- (-00.01,2);
\node at (0., -1.7) {0};
\node at (7,0.1) {\includegraphics[width=0.42\textwidth]{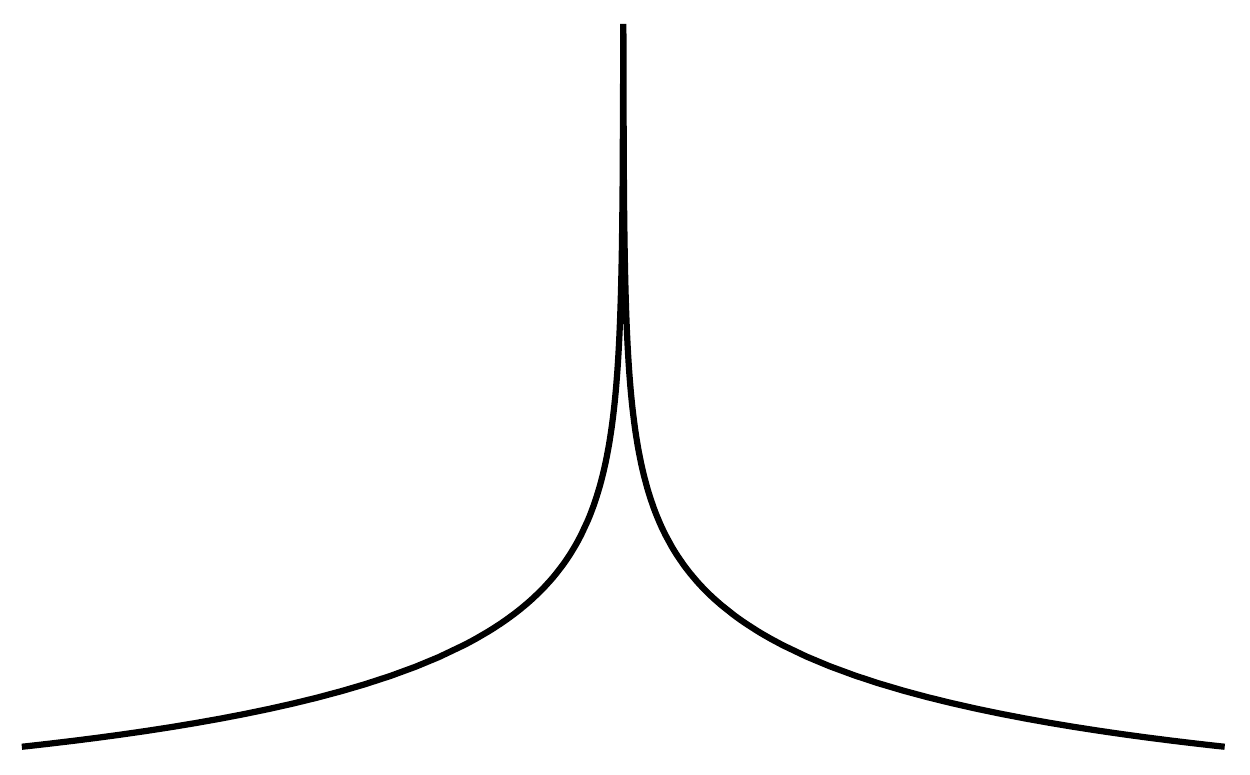}};
\draw [<->] (7-3,-1.52) -- (7+3,-1.52);
\draw [->] (6.99,-1.52) -- (6.99,2);
\node at (7, -1.7) {0};
\end{tikzpicture}
\caption{The radial profile of the convolution kernel $k_t(r)$ in $d=1$ dimensions (left) and $d=2$ (right) that we derive below.}
\end{figure}
\end{center}

We quickly comment on the distinguished role that the kernel plays. As indicated above, what is special about this kernel is that the error term depends linearly on $t$ (which is not too surprising: the average is taken over a ball of radius $\sim \sqrt{t}$, the linear part cancels) \textit{with a constant depending only on $\phi$ and $\|V\|_{L^{\infty}}$}. It is not difficult to see that by changing the kernel a tiny bit, one usually expects a change that is dependent on the size of $\Delta V$ (which may not even be defined since $V \in L^{\infty}$).

\begin{thm}[\cite{steini3}]
Let $\Omega \subset \mathbb{R}^n$ be an open, bounded domain with smooth boundary, let $0 \leq V \in C(\overline{\Omega})$ be a continuous potential and let $\phi$ be a solution of
\begin{align*}
(-\Delta + V)\phi &= \lambda \phi~ \quad \mbox{in~}\Omega \\
 \phi&= 0 \qquad \mbox{on}~ \partial \Omega.
\end{align*}
Then, for any fixed $x \in \Omega$, as $t \rightarrow 0$, we have, for $k_t$ as above,
$$ -\Delta \phi(x) +  (V*k_t)(x) \phi(x)= \lambda \phi(x) + \mathcal{O}_{\phi, \|V\|_{L^{\infty}}}(t),$$
where the implicit constant depends \emph{only} on $\phi$ and $\|V\|_{L^{\infty}}$. 
\end{thm} 

It is also shown that the landscape function $(-\Delta + V)u=1$ has a natural connection to this universal kernel and exhibits a similar type of stability property
$$ -\Delta u(x) + (V*k_t)(x) u(x)  = 1 + \mathcal{O}_{u, \|V\|_{L^{\infty}}}(t).$$
 Other methods for the purpose of fast computation of the location of localized low-lying eigenfunctions have been proposed \cite{alt, alt2, altmann, korn, lu, steini}. We also mention a curious localization phenomenon for Neumann boundary conditions \cite{sapo1, jones,  sapo2} that does seem to be of a different overall flavor.

\section{Fast Prediction of Localized Eigenstates}
We can now present the main contributions of this paper:
\begin{enumerate}
\item an alternative formulation of the universal convolution kernel that reduces the problem to two applications of the Fast Fourier Transform. This has a dramatic effect on computational
cost, allowing computations on a discrete $n \times n$ grid to be done in $\mathcal{O}(n^2 \log{n})$ allowing for much finer resolution.
\item exhaustive tests showing that $k_t * V$ has remarkable agreement with the regularized kernel $1/u$ derived from the landscape function and has comparable predictive power (partially explained by the results mentioned above).
\item and some results showing how the underlying idea can be extended
  to other operators, with special emphasis on the fractional
  Laplacian $(-\Delta)^{\alpha}$ and the bi-Laplacian $(-\Delta)^2$.
\end{enumerate}

We first discuss, in \S 2.1, the main idea behind the universal convolution kernel and how we will use it. \S 2.2 contains a precise description of the fast algorithm. \S 2.3 discusses the extension to other operators with special emphasis on the fractional Laplacian $(-\Delta)^{\alpha}$. Numerical Results are given in \S 3.

\subsection{The Main Idea} We first explain the main idea behind the universal kernel: in \cite{steini3}, this was derived using the Feynman-Kac formula, here we present a simpler motivation that only relies on the Duhamel principle. Let us suppose that
$$ (-\Delta + V) \phi = \lambda \phi.$$
We can then solve the heat propagator associated to the equation
$$ \frac{\partial u}{\partial t} - \Delta u=   -V u$$
starting with the initial condition $u(0,x) = \phi(x)$. 
Since eigenfunctions diagonalize that operator, we have $u(t,x) = e^{-\lambda t} \phi(x)$. In particular, the solution changes slowly and $u(t,x) \sim (1-\lambda t) \phi(x)$ for $t$ small.
At the same time, we can interpret this as an equation of the type $u_t - \Delta u = f$ and apply Duhamel's principle leading to
\begin{align*}
 u(t,x) &= e^{t\Delta} \phi + \int_0^t e^{(t-s)\Delta} f(s,x) ds.
\end{align*}
However, since we know the explicit solution, we have 
$$ f(s,x) = -V(x) u(s,x) = -V(x) e^{-\lambda s} \phi(x).$$
  We see that $\phi$ appears on all sides of the equation whereas $V$ appears in exactly one spot; what turns out to be relevant is the quantity
  $$ V_t =  \frac{1}{t} \int_{0}^{t} e^{(t-s) \Delta} V ds.$$
  However, $e^{t \Delta} f$ is simply the convolution with a Gaussian and thus $V_t$ can also be written as a convolution of averaged Gaussians. We refer to \cite{steini3} for details.

  \subsection{The Choice of $t$} In summary, we have described a natural kernel $k_t$ such that $V_t = k_t * V$ serves as a canonical mollification of the potential at scale $\sim \sqrt{t}$. It remains to understand how one should choose $t$. We emphasize that one of the advantages of our method being so fast is that one would reasonably check various values of $t$. It is clear from the motivation that we would like $\sqrt{t}$ (the scale of the convolution kernel $k_t$) to roughly correspond to the scale where most of the mass of the eigenfunction is supported -- the goal is to understand that scale without computing the eigenfunction first. Let us now suppose that $V \geq 0$ is given and that we have an eigenfunction
  $$ -\Delta \phi + V \phi = \lambda \phi.$$
  We can assume w.l.o.g. that the eigenfunction is mostly supported at scale $\sim r$ and is normalized in $L^2$, then multiplying by $\phi$ and integrating by parts shows
  $$ \int_{B(x_0,r)} |\nabla \phi|^2dx + \int_{B(x_0, r)} V(x) \phi(x)^2 dx = \lambda + \mbox{very small error},$$
  where the very small error comes from the mass of the eigenfunction outside the ball (which for the heuristic reasoning in this section we only need to be smaller than, say, $ \lambda/2$). Since these eigenfunctions are usually exponentially decaying, this is not a strong assumption.
  This computation does assume Dirichlet boundary conditions on $\partial \Omega$ but remains valid or very nearly valid in the general case for the same reason.
  We now try to understand the scaling of the first term: if $\phi$ is a bump function with $\| \phi\|_{L^{\infty}} = h$, then the $L^2-$normalization implies
  $$ h^2 r^d \sim 1 \qquad \mbox{and thus} \quad \int_{B(x_0,r)} |\nabla \phi|^2 dx\sim \left( \frac{h^2}{r^2} \right) r^d = \frac{1}{r^2}.$$
  At the same time, we expect that
  $$ \int_{B(x_0, r)} V(x) \phi(x)^2 dx \sim h^2\int_{B(x_0, r)} V(x) dx$$
  and we also expect that both terms should yield roughly the same contribution; if the gradient term was much larger, we would find a function with smaller energy by spreading out (note that in our setting, orthogonality tends to be less of an issue since we are only looking at localized low-frequency eigenfunctions). If the potential term was much larger, we may prefer to localize in a slightly smaller area (and preferably one where the potential is slightly smaller on average). Thus, we expect that
  $$ \frac{1}{r^2} \sim h^2 \int_{B(x_0, r)} V(x) dx.$$
  Multiplying both sides of the equation with $r^2$ and using $h^2 r^d \sim 1$ gives
  $$ 1 \sim \frac{1}{r^{d-2}}  \int_{B(x_0, r)} V(x) dx$$
  and thus suggests a natural scale definition of an eigenfunction at $x_0$ as
  $$ r(x_0, V) = \sup_{} \left\{ r>0:  \frac{1}{r^{d-2}}  \int_{B(x_0, r)} V(x) dx \leq 1\right\}.$$
  This is equivalent to a concept introduced by Fefferman \& Phong \cite{feff}. In particular, we expect that
  $$ \frac{1}{r(x_0, V)^2} \sim  \int_{B(x_0,r)} |\nabla \phi|^2 dx \sim \lambda.$$
  This has been made precise \cite{feff}. From a proof given by Shen \cite{shen}, we have
  $$ \int_{\Omega} \frac{|\phi(x)|^2}{r(x, V)^2} dx \lesssim_d  \int_{\Omega} |\nabla \phi|^2dx + \int_{\Omega} V(x) \phi(x)^2 dx.$$
  In particular, this tells us that when looking for low-frequency localized eigenfunction near a point $x_0$, we should select $t \sim r(x_0, V)^2$.
  
\subsection{The Fast Algorithm} We will now describe the fast algorithm. Since all the computations ultimately happen on a discretization of the underlying domain (see, for example, \cite{altmann} or \cite{arnold2}), we will describe the algorithm on a discretization of $\Omega = [0, 1)^2$ with periodic boundary conditions. The extension to higher dimensions will be obvious. We quickly note that this is, in some sense, the most general domain: the underlying impact of the potential to the localized eigenfunctions is ultimately local. It therefore does not matter whether we solve the problem on a square or on any other domain (and, for the same reason, it does not matter which boundary condition is imposed). In particular, we choose periodic boundary condition for the simplicity in applying fast Fourier transform (FFT).  Other boundary conditions can also be used with corresponding changes to the algorithm. For example, a fast discrete sine transform can be used for Dirichlet boundary conditions. \\

\textbf{Description of the Algorithm.} We assume the domain $[0,1)^2$ is
discretized with uniform grid with mesh size $h = \frac{1}{n}$, the resulting grid will be denoted by $\Omega_h$. We consider the discrete Laplacian
$\Delta_h$ defined by 
\begin{equation}
  (\Delta_h f)(x, y) = \frac{1}{h^2} \Bigl( f(x+h, y) + f(x-h, y) + f(x, y+h) + f(x, y-h) - 4f(x, y) \Bigr). 
\end{equation}
Define the discrete Fourier transform for
$(\xi, \eta) \in \Omega_h^{\ast} = (2\pi) \bigl\{0, 1, \ldots, n-1\bigr\}^2$
\begin{equation}
  \widehat{f}(\xi, \eta) = h^2 \sum_{(x, y) \in \Omega_h} e^{-i (\xi x + \eta y)} f(x, y). 
\end{equation}
The inverse Fourier transform is given by
\begin{equation}
  f(x, y) = \sum_{(\xi, \eta) \in \Omega_h^{\ast}} e^{ i (\xi x + \eta y)} \widehat{f}(\xi, \eta). 
\end{equation}
In particular, we observe that
\begin{equation}
  \begin{aligned}
    \widehat{\Delta_h f}(\xi, \eta) & = \frac{1}{h^2} \Bigl( e^{-i \xi h} +
    e^{i \xi h} + e^{- i \eta h} + e^{i \eta h} - 4 \Bigr) \widehat{f}(\xi,
    \eta) \\
    & = - \frac{4}{h^2} \Bigl( \sin^2\bigl(\frac{\xi h}{2}\bigr) +
    \sin^2\bigl(\frac{\eta h}{2}\bigr) \Bigr) \widehat{f}(\xi, \eta) =:
    M_h(\xi, \eta) \widehat{f}(\xi, \eta), 
  \end{aligned}
\end{equation}
where the last equality defines the Fourier multiplier associated with
the discrete Laplacian.
To construct the convolutional kernel corresponding to the discrete
operator, we consider the parabolic equation
\begin{align}\label{eq:parabolic}
  \partial_t u & = \Delta_h u,   &&\text{in } (0, \infty) \times \Omega_h; \\
  u(t = 0, \cdot) & = g, &&\text{on } \Omega_h. 
\end{align}
To obtain the corresponding semi-group, by taking the Fourier transform we have
\begin{equation}
  \partial_t \widehat{u}(t, \xi, \eta) = M_h(\xi, \eta) \widehat{u}(t, \xi, \eta), 
\end{equation}
and thus
\begin{equation}
  \widehat{u}(t, \xi, \eta) = e^{t M_h(\xi, \eta)} \widehat{u}(0, \xi, \eta) = e^{t M_h(\xi, \eta)} \widehat{g}(\xi, \eta). 
\end{equation}
The convolutional landscape funciton (at scale $t$) is given by
\begin{equation}\label{eq:WT}
  W_t(x, y) = h^2 \sum_{(x', y') \in \Omega_h} G_t(x - x', y - y') V(x', y'), 
\end{equation}
where the convolutional kernel $G_t$ is given by the time average of
the semigroup of the parabolic equation \eqref{eq:parabolic}. In the
Fourier space, it is thus 
\begin{equation}\label{eq:GT}
  \widehat{G}_t(\xi, \eta) = \frac{1}{t} \int_0^t e^{s M_h(\xi, \eta)} ds = \frac{1}{t M_h(\xi, \eta)} \Bigl( e^{t M_h(\xi, \eta)} - 1 \Bigr).
\end{equation}
The above is well defined for any non-zero
$(\xi, \eta) \in \Omega_h^{\ast}$; since we have periodic boundary
conditions, w.l.o.g. we can set $\widehat{G}_t(0, 0) = 1$.  The algorithm
for computing $W_t$ is as follows
\begin{itemize}
\item Given $V$ a grid function on $\Omega_h$, compute the discrete
  Fourier transform $\widehat{V}$ by FFT algorithm;
\item For each nonzero $(\xi, \eta) \in \Omega_h^{\ast}$, define
  \begin{equation}
    M_h(\xi, \eta) := - \frac{4}{h^2} \Bigl( \sin^2\bigl(\frac{\xi h}{2}\bigr) +
    \sin^2\bigl(\frac{\eta h}{2}\bigr) \Bigr),
  \end{equation}
  and calculate $\widehat{G}_t(\xi, \eta)$ according to \eqref{eq:GT};
\item For each nonzero $(\xi, \eta) \in \Omega_h^{\ast}$, calculate
  \begin{equation}
    \widehat{W}_t(\xi, \eta) = \widehat{G}_t(\xi, \eta) \widehat{V}(\xi, \eta);
  \end{equation}
\item Obtain $W_t$ by an inverse discrete Fourier transform applying on
  $\widehat{W}_t$ using FFT algorithm.
\end{itemize}

The main computational cost of the algorithm amounts to two
applications of FFT, and hence is of size $\mathcal{O}(n^2\log n)$.

\subsection{The Fractional Laplacian} Our approach is more widely applicable: the key ingredients are a good understanding of the short-time behavior 
of the associated parabolic equation, and such estimates are widely
available. We illustrate this with the fractional Laplacian
$$ (-\Delta)^{\alpha} u+ Vu = \lambda u.$$
Since we are working with the uniform grid in domain $[0, 1)^2$ with periodic boundary condition, we can define the fractional Laplacian spectrally to avoid the usual difficulties in properly defining the fractional Laplacian on a bounded
domain. For any $\alpha \geq 0$, we define the fractional Laplacian as the operator that sends
$$ (-\Delta)^{\alpha} e^{ i(\xi x + \eta y)} = \bigl( \lvert\xi\rvert^2 + \lvert\eta\rvert^2 \bigr)^{\alpha}  e^{ i(\xi x + \eta
  y)}$$
for complex exponentials with
$(\xi, \eta) \in \Omega_h^{\ast} = (2\pi) \bigl\{-N/2, -N/2 + 1,
\ldots, N/2-1 \bigr\}^2$.  Naturally, we recover the standard
Laplacian (with pseudospectral discretization) for $\alpha = 1$. We see that $0 < \alpha < 1$ impacts the
extent to which frequencies get dampened: the smaller $\alpha$, the
more contribution we get from small frequencies. Formulated on the
spatial side, the operator $(-\Delta)^{\alpha}$ is non-local and
$\alpha$ governs the scale of the of non-locality. The smaller
$\alpha$, the stronger the non-locality and the harder it becomes for
eigenstates to localize.  Nonetheless, our method remains
applicable. Setting
$$ M_h(\xi, \eta) = - \bigl( \lvert\xi\rvert^2 + \lvert\eta\rvert^2 \bigr)^{\alpha},$$
we have the corresponding multiplier completely unchanged
\begin{equation}\label{eq:GT2}
  \widehat{G}_t(\xi, \eta) = \frac{1}{t} \int_0^t e^{s M_h(\xi, \eta)} ds = \frac{1}{t M_h(\xi, \eta)} \Bigl( e^{t M_h(\xi, \eta)} - 1 \Bigr)
\end{equation}
and we can proceed as above. We refer to \S\ref{sec:fraclap} for
numerical examples.

\section{Numerical Results}
This section discusses a variety of numerical results. We first show, in \S 3.1, that the regularized potential $W_t = k_t * V$ approximates the landscape function fairly well. This is further evidence
that the two notions are, on some level, connected (see also \cite[Theorem 2]{steini3}). In particular, if one were to solve $(-\Delta + V)u = 1$ by means of an iterative method, $W_t$ may be
a good initial value to choose from. In \S 3.2. we investigate to what extent the local minima of $W_t$ do indeed correspond to low-frequency eigenfunctions and, conversely, analyze how many low-frequency eigenfunctions are being detected by the local minima.

\subsection{Setup} We use the same uniform setup for all examples: computations were carried out on a $256 \times 256$ grid with boundary identified to obtain the Torus geometry. The potential $V$ at each point is an i.i.d. random variable chosen uniformly from $[0, V_{\max}]$.  For small values of $V_{\max}$, the corresponding eigenfunctions are quite delocalized -- as $V_{\max}$ increases, eigenfunctions become more localized and we have chosen parameters to cover that entire range, this happens to be $1 \leq V_{\max} \leq 5$. 
For each setting of $V_{\max}$, we simulated 100 random instances on which we solved the eigenvalue problem, computed the Filoche-Mayboroda Landscape function $1/u$ and our regularized potential $W_t$ for various choices of $t$. In each instance, we compute the first 64 true eigenfunctions for $-\Delta + V$. For both the landscape function as well as the
regularized potential, we look for the 16 smallest local minima and ordered them increasing in size.  We say that a local minimum has found an eigenfunction whenever the grid point in which the minimum occurs and the grid point in which an eigenfunction has its largest absolute value are at most $5h$ (in Euclidean distance) apart (recall that $h$ is the mesh size); we then associate that eigenfunction to
that minimum. This will usually lead to us localizing 16 eigenfunctions (which, by default, are always taken from the first 64).\\

We used the following test-statistics to quantify the performance of various approaches:
\begin{enumerate}
\item Eigenvalue Ratios (\textsc{EigRat})
\item The First Missing Eigenfunction (\textsc{FirstMissEig})
\item The First Unused Minimum (\textsc{FirstMissMin})
\item The Number of Unused Minima (\textsc{DismissedMin})
\end{enumerate}

\subsubsection{The Eigenvalue Ratio.} Given our 16 local minima, we
associate a number of eigenfunctions $\phi_{i_1}, \dots, \phi_{i_k}$ to them.  We would ideally like to have
$k=16$, but it certainly occurs that some of the first minima do not
correspond to any of the first $64$ eigenfunctions; in that case $k$ is the number of minima with corresponding eigenvalues. Moreover, we would
ideally like $i_1=1, i_2 = 2, ... $. However, this again is a
misleading metric since the quantity relevant to the underlying
physical system is, in actuality, given by the eigenvalue. Thus, for each such
realization, we define the eigenvalue ratio \textsc{EigRat} via 
$$ \mbox{EigRat} = \frac{\lambda_{i_1} + \dots + \lambda_{i_k}}{\lambda_1 + \dots + \lambda_k}.$$
This number is always at least $1$. The smaller this number is, the closer it is to $1$, the better we perform at finding \textit{low-frequency} localized eigenstates.

\subsubsection{The First Missing Eigenfunction.} Once we have found our 16 local minima and associated (at most) 16 eigenfunctions to it. However, these eigenfunctions will, generally, not be in order and there is certainly an eigenfunction $\phi_k$ which is \textit{not} being found by any of the 16 local minima. \textsc{FirstMissEig} simply returns the smallest value $k$ for which this happens, the smallest eigenfunction that is being missed by the first 16 minima. The larger \textsc{FirstMissEig}, the better the method performs at finding all small eigenfunctions. 

\subsubsection{The First Unused Minimum.} After having found the first $16$ minima and ordering them in increasing size, it will usually happen that some minimum will not correspond to any nearby eigenfunction. \textsc{FirstMissMin} lists the smallest (as ordered by size) index of a minimum that does not correspond to one of the first $64$ eigenfunctions.  The smaller the value in the minimum, the more strongly we expect it to correspond to a small eigenfunction and thus, for \textsc{FirstMissMin}, the larger the value, the better.

\subsubsection{The Number of Unused Minima.} In a similar spirit as \textsc{FirstMissMin}, \textsc{DismissedMin} counts the number of the first 16 minima that do not end up corresponding to one of the first $64$ eigenfunctions. The smaller the number, the better.\\

\begin{figure}[h!]
\begin{minipage}[l]{.49\textwidth}
\begin{tikzpicture}
\node at (0,0) {\includegraphics[width = 1\textwidth]{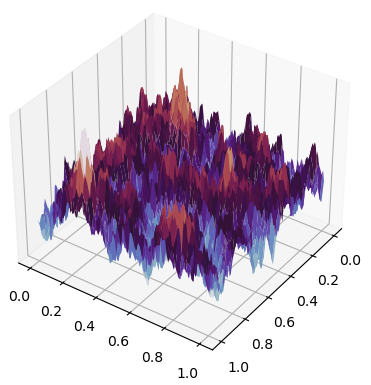}};
\end{tikzpicture}
\end{minipage} 
\begin{minipage}[r]{.49\textwidth}
\begin{tikzpicture}
\node at (0,0) {\includegraphics[width = 1\textwidth]{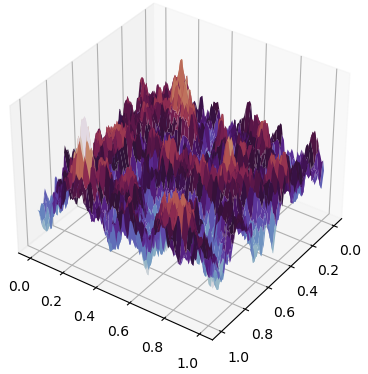}};
\end{tikzpicture}
\end{minipage} 
\caption{The landscape function and $W_t$. We see that the resulting pictures are visually quite similar \label{fig:comp} (see \S 3.2).} 
\end{figure}

These statistics are designed to capture aspects of an ideal
approximating potential: it finds most low-frequency eigenfunctions
accurately and roughly in the right order (where errors are to be
measured in terms of the actual eigenvalues as opposed to the ordering
itself). Moreover, most minima are supposed to correspond to local
eigenfunctions. We emphasize that there are many
other notions of quality that could be used.

\begin{table}[h!]
\begin{center}
\begin{tabular}{  c c c c c } 
 $V_{\max}=2$  & $\| \cdot - \cdot\|_{L^1}$ & $\|\cdot - \cdot\|_{L^2}$ & $\| \cdot - \cdot \|_{L^{\infty}} $ \\ 
   \hline
t=0.10 & 0.092 &  0.120 &  0.446  \\
t=0.25 &  0.071 &  0.087 & 0.318  \\
t=0.40 & 0.065 & 0.075  & 0.233   \\
\hline
\hline
 $V_{\max}=5$  & $\| \cdot - \cdot\|_{L^1}$ & $\|\cdot - \cdot\|_{L^2}$ & $\| \cdot - \cdot \|_{L^{\infty}} $ \\ 
   \hline
t=0.5 & 0.067 &  0.083 &  0.309  \\
t=1.0 &  0.058 &  0.071 & 0.252  \\
t=1.5 & 0.053 & 0.063  & 0.211   \\
\hline
\end{tabular}
\end{center}
\caption{The average difference between (normalized) landscape function and (normalized) regularized potential $W_t$ (see \S 3.2).}
\end{table}

\vspace{-20pt}

\subsection{Approximation to $1/u$}
The purpose of this section is to point out that, visually, both the landscape function as well as our regularization are often quite similar (see Fig.~\ref{fig:comp}). The effect is rather striking in $d=1$ dimensions (see the examples in \cite{steini3}).
We briefly compare this by considering the $L^p$ norm between the two types of regularized potentials.  Since we eventually end up localizing eigenfunctions
by looking at the location of minima of the respective functions, we normalize them to assume minimal value 0 and maximal value 1.  A comparison between
the two functions shows that they are typically quite close.
We see that the similarity in all norms increases with larger values of $t$. 
One should interpret this with great care: for the purpose of localization, we care about the location and the value of the minima -- this is not easily captured by
comparing $L^p-$norms. However, the results discussed in this section do suggest that both methods, the landscape function and the regularized
potential, capture the same type of object.

\subsection{Summary of Results} We observe that our approach yields satisfying results across a wide range of parameters; indeed, we observe that, in the setting we study, our approach generally leads to results that are superior to that of the landscape functions in the setting where eigenfunctions are localized over a wider area $(V_{\max} \in \left\{1, 1.5, 2\right\}$). For more localized eigenfunctions $(V_{\max} \in \left\{5\right\})$, the landscape function produces more accurate results. Both methods seem to be more or less comparable -- both tend to find low-frequency eigenfunctions with remarkable efficiency (the values of $\textsc{EigRat}$ are very close to 1). Our method has also proven remarkably robust for many values of $t$ (there is an overall tendency that picking the scale of the kernel to be roughly at the scale of the eigenfunction leads to slightly better results, we refer to \S 2.2. for what we would expect this scale to be in the setting of a random potential).

\subsubsection{Fairly delocalized eigenfunctions, $V_{\max} = 1$.} We start with the case where each lattice point assumes a random value in $[0,1]$. We see (in Fig. 3) that this leads to dispersed eigenfunctions that are spread out over a large area.
\begin{figure}[h!]
\begin{minipage}[l]{.495\textwidth}
\begin{tikzpicture}
\node at (0,0) {\includegraphics[width = 0.9\textwidth]{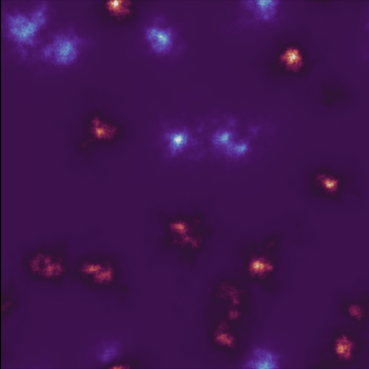}};
\end{tikzpicture}
\end{minipage} 
\begin{minipage}[r]{.495\textwidth}
\begin{tikzpicture}
\node at (0,0) {\includegraphics[width = 0.9\textwidth]{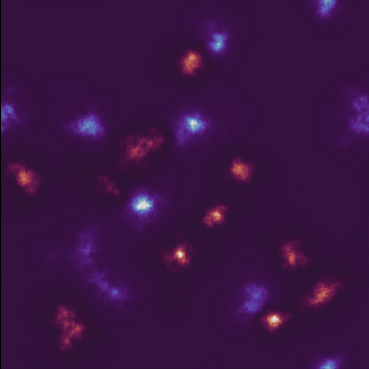}};
\end{tikzpicture}
\end{minipage} 
\caption{Two instances of $V_{\max} =1$ and the first 20 eigenfunctions on each (colored by sign). } 
\end{figure}
\begin{table}[h!]
\begin{center}
\begin{tabular}{  c c c c c } 
   & \textsc{EigRat} & \textsc{FirstMissEig} & \textsc{FirstMissMin} & \textsc{DismissedMin} \\ 
   \hline
 Landscape $1/u$ & 1.0094 & 3.19 & 6.55 & 3.32 \\ 
 \hline
t=0.01 & 1.0103 & 2.13 & 1.66 & 3.32 \\
t=0.05 & 1.0030 & 7.35 & 6.39 & 2.84 \\
t=0.10 & 1.0031 & 8.01 & 8.20 & 1.49 \\
t=0.20 & 1.0063 & 5.45 & 8.80 & 1.80 \\
\hline
\end{tabular}
\end{center}
\caption{Statistics for $V_{\max} = 1$.}
\end{table}
We see that, for all values of $t$, our method finds eigenfunctions corresponding to small eigenvalues -- however, it rarely finds all of them and $\textsc{FirstMissEig}$, the average index of the first eigenfunction not uncovered by the first 16 minima, is quite small. However, the moment $t$ corresponds to slightly more localized kernels, we obtain very strong results: the first 16 local minima uncover, on average, the first 7 eigenfunctions exactly and, moreover, the first 6 minima do correspond to one of the first few eigenfunctions.

\subsubsection{Moderately delocalized eigenfunctions, $V_{\max} = 1.5$.} 
Our next case deals with eigenfunctions that may still be spread out over larger regions, but tend to be more concentrated on average (see Fig. 4). We observe again that both the landscape as well as our method perform remarkably well. For a suitable range of parameters, virtually all minima correspond to local eigenfunctions. In that regime we also observe that our method finds twice as many low-frequency eigenfunctions as the landscape function. \textsc{EigRat} is uniformly small.
\begin{figure}[h!]
\begin{minipage}[l]{.495\textwidth}
\begin{tikzpicture}
\node at (0,0) {\includegraphics[width = 0.9\textwidth]{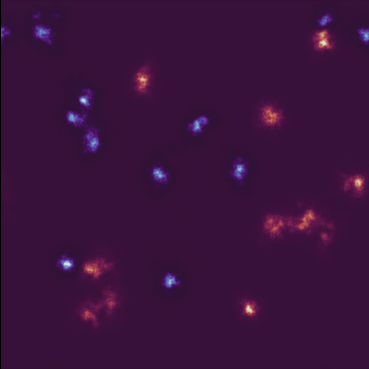}};
\end{tikzpicture}
\end{minipage} 
\begin{minipage}[r]{.495\textwidth}
\begin{tikzpicture}
\node at (0,0) {\includegraphics[width = 0.9\textwidth]{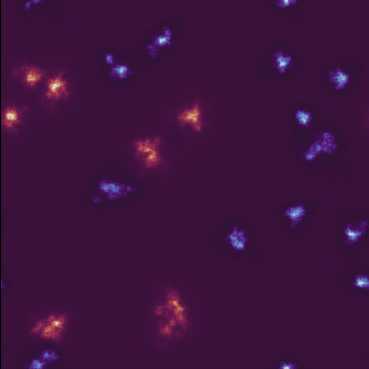}};
\end{tikzpicture}
\end{minipage} 
\caption{Two instances of $V_{\max} =1.5$ and the first 20 eigenfunctions on each (colored by sign).} 
\end{figure}

\begin{table}[h!]
\begin{center}
\begin{tabular}{  c c c c c } 
   & \textsc{EigRat} & \textsc{FirstMissEig} & \textsc{FirstMissMin} & \textsc{DismissedMin} \\ 
   \hline
 Landscape $1/u$ & 1.0097 & 4.17 & 8.01 & 2.52 \\ 
 \hline
t=0.06 & 1.0071 & 4.50 & 6.30 & 2.42 \\
t=0.12 & 1.0032 & 7.91 & 8.27 & 1.03 \\
t=0.18 & 1.0036 & 7.9 & 9.17 & 0.81 \\
t=0.24 & 1.0050 & 6.92 & 10.05 & 0.87 \\
t=0.30 & 1.0061 & 6.13 & 10.39 & 1.11\\
\hline
\end{tabular}
\end{center}
\caption{Statistics for $V_{\max} = 1.5$.}
\end{table}

\subsubsection{Localized eigenfunctions, $V_{\max} = 2$.} 
The next case, the potential assuming i.i.d. values uniformly in $[0,2]$, corresponds to even better behaved eigenfunctions which now exhibit a strong-form concentration (see Fig. 5): very few are spread out, though all of them extend over a certain nontrivial region. Again, we observe stability and efficiency of our method for a wide range of parameters.
\begin{figure}[h!]
\begin{minipage}[l]{.495\textwidth}
\begin{tikzpicture}
\node at (0,0) {\includegraphics[width = 0.9\textwidth]{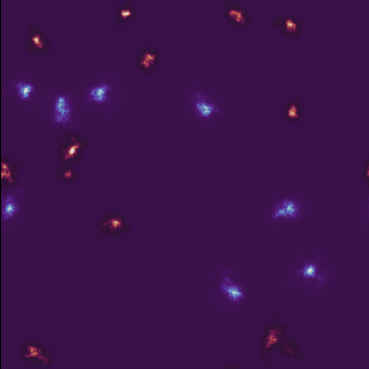}};
\end{tikzpicture}
\end{minipage} 
\begin{minipage}[r]{.495\textwidth}
\begin{tikzpicture}
\node at (0,0) {\includegraphics[width = 0.9\textwidth]{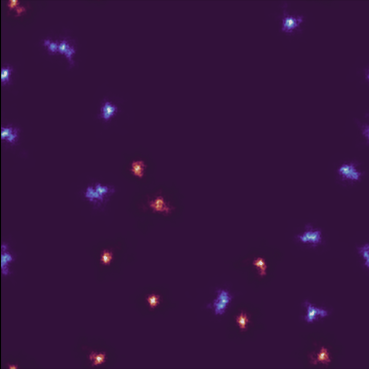}};
\end{tikzpicture}
\end{minipage} 
\caption{Two instances of $V_{\max} =2$ and the first 20 eigenfunctions on each (colored by sign). } 
\end{figure}

\begin{table}[h!]
\begin{center}
\begin{tabular}{  c c c c c } 
   & \textsc{EigRat} & \textsc{FirstMissEig} & \textsc{FirstMissMin} & \textsc{DismissedMin} \\ 
   \hline
 Landscape $1/u$ & 1.0097 & 5.13 & 9.64 & 1.80 \\ 
 \hline
t=0.1 & 1.0077 & 4.77 & 8.76 & 1.57 \\
t=0.2 & 1.0045 & 7.98 & 9.84 & 0.67 \\
t=0.3 & 1.0052 & 7.68 & 10.76 & 0.49 \\
t=0.4 & 1.0068 & 6.76 & 11.06 & 0.73 \\
\hline
\end{tabular}
\end{center}
\caption{Statistics for $V_{\max} = 2$.}
\end{table}

\subsubsection{Highly Localized eigenfunctions, $V_{\max} = 5$.} 
Our next case corresponds to a much larger potential, with values assigned uniformly from $[0,5]$, and this comes with more concentrated eigenfunctions that are now strongly localized (many of them carry most of their mass on a $8 \times 8$ grid). We observe that the landscape function produces slightly more accurate results. Both methods show remarkably large values of \textsc{FirstMissMin}, which means that the minima do correspond quite strongly to the ordering of the eigenvalues. This is, in some sense, the easiest case: the potential is quite large which leads to strong localization properties -- in particular, the only relevant information whether there might be a localized eigenfunction near any given point is the behavior of $V$ in a small neighborhood (`small' compared to what it would be for smaller potentials).
\begin{center}
\begin{figure}[h!]
\begin{minipage}[l]{.495\textwidth}
\begin{tikzpicture}
\node at (0,0) {\includegraphics[width = 0.9\textwidth]{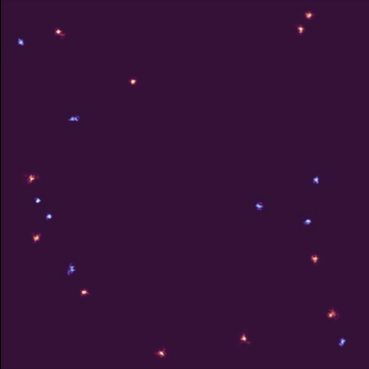}};
\end{tikzpicture}
\end{minipage} 
\begin{minipage}[r]{.495\textwidth}
\begin{tikzpicture}
\node at (0,0) {\includegraphics[width = 0.9\textwidth]{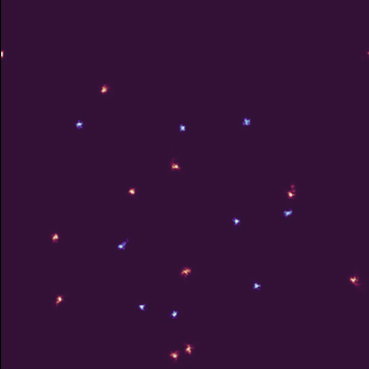}};
\end{tikzpicture}
\end{minipage} 
\caption{Two instances of $V_{\max} =5$ and the first 20 eigenfunctions on each (colored by sign).} 
\end{figure}
\end{center}

\begin{table}[h!]
\begin{center}
\begin{tabular}{  c c c c c } 
   & \textsc{EigRat} & \textsc{FirstMissEig} & \textsc{FirstMissMin} & \textsc{DismissedMin} \\ 
   \hline
 Landscape $1/u$ & 1.0077 & 8.16 & 12.12 & 0.43 \\ 
 \hline
t=0.50 & 1.0149 & 4.60 & 9.89 & 1.35 \\
t=0.75 & 1.0125 & 6.06 & 11.05 & 0.93 \\
t=1.00 & 1.0122 & 6.20 & 10.91 & 0.81 \\
t=1.25 & 1.0121 & 6.18 & 11.16 & 1.07 \\
\hline
\end{tabular}
\end{center}
\caption{Statistics for $V_{\max} = 5$.}
\end{table}

\subsection{Comparison with Gaussian Filtering} 
 Our approach is based on the idea that taking suitable local averages of the potential can create a reasonable
enough approximation of the effective potential so that one can read off the location of localized low-frequency eigenstates from the minimum of this convolution.
Our particular choice of convolution kernel is motivated by theory, however, it does raise the question of what happens if we merely take a Gaussian (note also
that our convolution kernel is a specific superposition of Gaussians). We observe that it suffices to change Equation (\ref{eq:GT}) in the algorithm to
\begin{equation}
  \widehat{G}_t(\xi, \eta) =  e^{t M_h(\xi, \eta)}. 
\end{equation}

We observe that simple Gaussian filtering by itself is also remarkably effective; perhaps this is not all that surprising considering that the Gaussian is actually the short time asymptotic profile
of the associated parabolic equation, and seeing as well that our kernel $k_t$ is actually a superposition of Gaussians. 

\begin{table}[h!]
\begin{center}
\begin{tabular}{  c c c c c } 
 $V_{\max}=1$  & \textsc{EigRat} & \textsc{FirstMissEig} & \textsc{FirstMissMin} & \textsc{DismissedMin} \\ 
   \hline
 Landscape $1/u$ & 1.0094 & 3.61 & 6.6 & 3.18 \\ 
 \hline
$k_{0.1}*V$ & 1.0031 & 8.01 & 8.20 & 1.49 \\
\hline
$g_{0.2}*V$ & 1.0035 & 6.55 & 7.51 & 1.94 \\
$g_{0.3}*V$ & 1.0040 & 6.74 & 8.87 & 1.30 \\
$g_{0.4}*V$ & 1.0055 & 5.56 & 9.31 & 1.38 \\
\hline\\
 $V_{\max}=2$  & \textsc{EigRat} & \textsc{FirstMissEig} & \textsc{FirstMissMin} & \textsc{DismissedMin} \\ 
   \hline
 Landscape $1/u$ & 1.0096 & 5.27 & 10.14 & 1.64 \\ 
 \hline
$k_{0.2}*V$ & 1.0045 & 7.98 & 9.84 & 0.67 \\\hline
$g_{0.3}*V$ & 1.0092 & 4.56 & 9.23 & 1.52 \\
$g_{0.6}*V$ & 1.0048 & 7.58 & 9.67 & 0.75 \\
$g_{0.9}*V$ & 1.0061 & 6.65 & 10.34 & 0.65 \\
\hline\\
 $V_{\max}=5$  & \textsc{EigRat} & \textsc{FirstMissEig} & \textsc{FirstMissMin} & \textsc{DismissedMin} \\ 
   \hline
 Landscape $1/u$ & 1.0096 & 5.27 & 10.14 & 1.64 \\ 
 \hline
$k_1*V$ & 1.0122 & 6.20 & 10.91 & 0.81 \\
\hline
$g_{1}*V$ & 1.0197 & 3.4 & 9.38 & 2.11 \\
$g_{2}*V$ & 1.0116 & 5.91 & 11.07 & 0.87 \\
$g_{3}*V$ & 1.0141 & 4.97 & 10.55 & 1.41 \\
\end{tabular}
\end{center}
\caption{Statistics for the landscape, the convolution with our kernel $k_t$ and the convolution with a Gaussian $g_t$.}
\end{table}
Nonetheless, we do observe consistently that our kernel is better at localizing low-frequency eigenfunctions (\textsc{EigRat} being smaller and \textsc{FirstMissEig} being larger). However, when it comes to having minima roughly correspond with some low-frequency eigenfunctions (\textsc{DismissedMin}) the Gaussian sometimes yields the best results.
\begin{table}[h!]
\begin{center}
\begin{tabular}{  c c c c c } 
 $V_{\max}=2$  & \textsc{EigRat} & \textsc{FirstMissEig} & \textsc{FirstMissMin} & \textsc{DismissedMin} \\ 
   \hline
 Landscape $1/u$ & 1.0094 & 3.61 & 6.6 & 3.18 \\
 \hline
$k_{0.2}*V$ & 1.0045 & 7.98 & 9.84 & 0.67 \\
$g_{0.6}*V$ & 1.0048 & 7.58 & 9.67 & 0.75 \\
\hline
$3 \times 3$ box & 1.0223 & 2.03 & 3.38 & 6.13 \\
$5 \times 5$ box & 1.0127 & 3.87 & 7.44 & 2.8  \end{tabular}
\end{center}
\caption{Statistics for the landscape, the convolution with our kernel $k_t$ and a Gaussian $g_t$ and compared to the approximation of the effective potential arising from averages over pixel boxes.}
\end{table}
The Gaussian kernel and our convolution kernel being so similar, we chose to compare with another method of averaging that is distinctly non-Gaussian in nature: computing the average in surrounding $3 \times 3$ and $5 \times 5$ pixel boxes and using these averages as an effective potential; seeing as it is the local behavior of $V$ that enables localization, this is not an unreasonable way of computing an average. We see that the method works, but that it behaves quite differently from our method.

\subsection{The Fractional Laplacian}\label{sec:fraclap}
As discussed in \S 2.3, our method
naturally extends to the fractional Laplacian as well. We will quickly
discuss examples for the operator $(-\Delta)^{3/4} + V$, where $V$ is
again an i.i.d. random potential assuming random values, uniformly in
$[0, V_{\max}]$. We will compare the performance of our method against
the \textit{fractional} Landscape $1/u$, where $u$ solves
$$ \left[(-\Delta)^{\alpha} + V\right]u = 1.$$

 The problem is much harder due to strong non-locality in the operator. We adjust $V_{\max}$ such that the
localization is roughly equivalent to the standard Laplacian cases
with $V_{\max} = 1, 2, 5$. To find this rough equivalence, we examine
the heat-maps generated from our standard Laplacian cases and choose
$V_{\max}$ such that the local peaks are roughly the same diameter for
$\alpha = 3/4$. The spreading of the eigenfunctions beyond
those peaks are much larger for the fractional Laplacian case. We
found comparable $V_{\max}$ values for the fractional Laplacian at
$\alpha = 3/4$ to be $3/40, 5/40$, and $10/40$.
Since the discretized fractional Laplacian leads to a dense matrix, we define it as an operator acting on vectors and use iterative solvers for the Landscape function and eigenfunctions. 

\begin{table}[h!]
\begin{center}
\begin{tabular}{  c c c c c } 
 $V_{\max}= 3/40$  & \textsc{EigRat} & \textsc{FirstMissEig} & \textsc{FirstMissMin} & \textsc{DismissedMin} \\ 
   \hline
  fractional Landscape & 1.014 & 2.02 & 4.66 & 5.24 \\
$t=0.01$ & 1.009 & 3.41 & 6.34 & 3.51 \\
$t=0.02$ & 1.011 & 2.55 & 6.05 & 4.13 \\
\hline
 $V_{\max}=5/40$  & \textsc{EigRat} & \textsc{FirstMissEig} & \textsc{FirstMissMin} & \textsc{DismissedMin} \\ 
   \hline
     fractional Landscape & 1.009 & 4.32 & 9.01 & 2.28 \\
$t=0.01$ & 1.005 & 6.2 & 7.72 & 1.83 \\
$t=0.02$ & 1.005 & 6.99 & 8.91 & 1.08 \\
\hline
 $V_{\max}=10/40$  & \textsc{EigRat} & \textsc{FirstMissEig} & \textsc{FirstMissMin} & \textsc{DismissedMin} \\ 
   \hline
        fractional Landscape & 1.009 & 6.73 & 11.7 & 0.61 \\
$t=0.02$ & 1.017 & 3.6 & 8.38 & 2.05 \\
$t=0.05$ & 1.011 & 5.65 & 11.2 & 1.13 \\
\end{tabular}
\end{center}
\caption{Statistics for the non-local problem $(-\Delta)^{3/4} + V$.}
\end{table}

We observe behavior that is consistent with the purely Laplacian case: as the potential increases, eigenfunctions
become increasingly localized and, therefore, easier to detect. Moreover, local minima in the regularized
potential carry more and more information. As in the Laplacian case, the fractional landscape becomes more accurate
as eigenfunctions become more localized.

\subsection{The Bi-Laplacian}

We also test the bi-Laplacian operator $(-\Delta)^2+V$, where $V$ is again an i.i.d. random potential assuming values uniformly at random in $[0, V_{\max}]$. We will again compare this to an adapted landscape function $1/u$, where
$$ \left[ (-\Delta)^2 + V \right] u = 1.$$

 Compared with the Laplacian or the fractional Laplacian examples tested, the bi-Laplacian operator is much more singular. We again adjust $V_{\max}$ such that the localization is roughly equivalent to the standard Laplacian trials with $V_{\max} = 1, 2, 5$. By similar procedure as for the fractional Laplacian case, we find comparable $V_{\max}$ values for the bi-Laplacian to be $1/200000$, $4/200000$, and $10/200000$. Note that the potential value is much smaller due to the stronger singularity of the bi-Laplacian operator. 
For discretization, we simply take the square of the finite difference approximation of the Laplacian operator, which is still quite sparse.

\begin{table}[h!]
\begin{center}
\begin{tabular}{  c c c c c } 
 $V_{\max}= 1/200000$  & \textsc{EigRat} & \textsc{FirstMissEig} & \textsc{FirstMissMin} & \textsc{DismissedMin} \\ 
   \hline
   Landscape $1/u$ & 1.011 & 3.73 & 5.97 & 3.42 \\
   $(-\Delta)^2-$ Landscape & 1.006 & 6.7 & 1.06 & 8.3 \\
$t=2 \cdot 10^{-3}$ & 1.005 & 5.85 & 7.73 & 1.65 \\
$t=6 \cdot 10^{-3}$ & 1.001 & 10.05 & 9.80 & 0.5 \\
\hline
 $V_{\max}=4/200000$  & \textsc{EigRat} & \textsc{FirstMissEig} & \textsc{FirstMissMin} & \textsc{DismissedMin} \\ 
   \hline
      Landscape $1/u$ & 1.022 & 2.7 & 5.09 & 5.06 \\
   $(-\Delta)^2-$ Landscape & 1.005 & 7.92 & 11.34 & 0.24 \\
$t=1 \cdot 10^{-2}$ & 1.013 & 3.94 & 10.28 & 1.02 \\
$t=1 \cdot 10^{-1}$ & 1.006 & 7.73 & 11.78 & 0.32 \\
\hline
 $V_{\max}=10/200000$  & \textsc{EigRat} & \textsc{FirstMissEig} & \textsc{FirstMissMin} & \textsc{DismissedMin} \\ 
   \hline
      Landscape $1/u$ & 1.037 & 1.83 & 3.6 & 5.83 \\
   $(-\Delta)^2-$ Landscape & 1.005 & 9.03 & 10.6 & 0.19 \\
$t=1 \cdot 10^{-1}$ & 1.008 & 7.54 & 10.90 & 0.22 \\
$t=4 \cdot 10^{-1}$ & 1.014 & 5.41 & 10.22 & 1.26 \\
\end{tabular}
\end{center}
\caption{Statistics for the non-local problem $(-\Delta)^2 + V$.}
\end{table}

We observe again results that are fairly consistent with what we have seen in other cases: as the potential
increases, eigenfunctions become more localized and the corresponding problem becomes easier. We see that
the classical landscape function does not do very well (nor is there any reason it should, it was not designed for
this problem). In contrast, the landscape function adapted to the Bi-Laplacian does indeed do quite well and becomes
better as the problem becomes more localized. Again, as we have seen many times, our method outperforms other
approaches in the case where eigenfunctions are fairly delocalized. \\

\begin{figure}[h!]
\begin{minipage}[l]{.49\textwidth}
\begin{tikzpicture}
\node at (0,0) {\includegraphics[width = 1\textwidth]{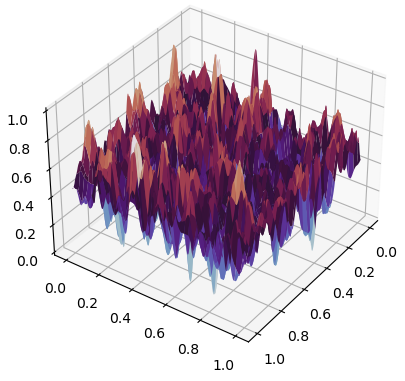}};
\end{tikzpicture}
\end{minipage} 
\begin{minipage}[r]{.49\textwidth}
\begin{tikzpicture}
\node at (0,0) {\includegraphics[width = 1\textwidth]{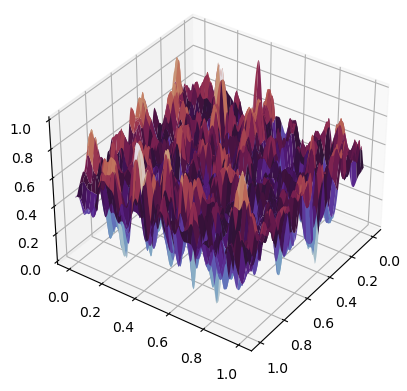}};
\end{tikzpicture}
\end{minipage}
\begin{minipage}[c]{.49\textwidth}
\begin{tikzpicture}
\node at (0,0) {\includegraphics[width = 1\textwidth]{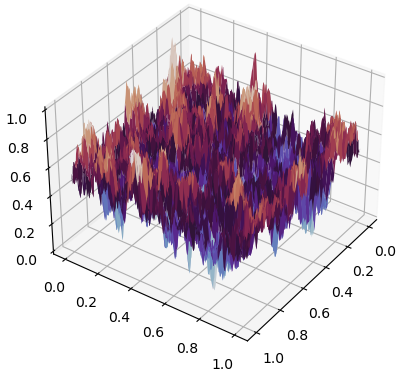}};
\end{tikzpicture}
\end{minipage} 
\caption{The bi-Laplacian landscape function (top left), $W_t$ for bi-Laplacian (top right), and landscape function (bottom). \label{fig:fraccomp}}
\end{figure}

In addition to using the landscape function corresponding to the
bi-Laplacian, given by $((-\Delta)^2 + V ) u = 1$, we also compute the
landscape function using the standard Laplacian (with the same
potential), see Fig.~\ref{fig:fraccomp} for comparison. We observe that the
bi-Laplacian landscape function and $W_t$ are qualitatively quite similar,
while the Laplacian landscape function is much more oscillatory, even though it still captures the location of large peaks.

\end{document}